\documentstyle{ioplppt}
\newcommand{\R}{{\sf R\hspace*{-0.9ex}\rule{0.15ex}%
    {1.5ex}\hspace*{0.9ex}}}
\begin{document}
\jl{1}

\def \Ha #1 { {\cal H}_{ #1 } }
\def \Ho { {\Ha 1 } }
\def \Ht { {\Ha 2 } }
\def \M {M^{3n}}
\def \th {{(3)}}
\def \tw {{(2)}}
\def \a {{\alpha}}
\def \b {{\beta}}
\def \omegath {{\omega^{\th}}}
\def \omegatw {{\omega^{\tw}}}
\def \div #1 #2 {{ { #1 } \over { #2 } }}
\def \dbd #1 #2 { { \div {\partial { #1 }} {\partial { #2 }} }}
\def \pf {{\bf Proof: }}
\def \ex {{\bf Example: }}
\def \T {{\cal T}}
\def \X {{\cal X}}
\def \S {{\cal S}}
\def \Xl {{\cal X}_L}
\def \Xr {{\cal X}_R}
\def \g {{\makebox{\boldmath $g$}}}
\def \L {{\cal L }}
\def \RR {{\cal R}}
\def \qed {{\rightline{$\Box$}}}

\newtheorem{defn}{Definition}[section]
\newtheorem{thm}{Theorem}[section]
\newtheorem{prop}{Proposition}[section]
\newtheorem{lemma}{Lemma}[section]
\newtheorem{cor}{Corollary}[thm]

\title{Momentum maps and Noether theorem for generalized Nambu mechanics} 
\author{ Sagar A. Pandit\dag\ftnote{4}{e-mail: sagar@prl.ernet.in} and 
  Anil D. Gangal\ddag\ftnote{5}{e-mail: adg@physics.unipune.ernet.in}}

\address{\dag Physical Research Laboratory, Navrangpura, Ahmedabad 380 009, India}
\address{\ddag Department of Physics, University of Pune, Pune 411 007, India.}
\begin{abstract}
In Ref.~\cite{Sag}  we proposed  a geometric formulation of
generalized Nambu mechanics. In the present paper we extend the class
of Nambu systems by replacing the stringent condition of constancy of
3-form by closedness. We also explore the connection between
continuous groups of symmetries and conservation laws for such
systems. The Noether theorem for generalized Nambu systems is
formulated by generalizing the notion of momentum map. In this case, a
natural choice of dynamical variables for discussion of symmetries is
2-form fields. Thus the generators and the conserved quantities in
Noether theorem are best expressed in terms of 2-forms. The connection
between the generators and the conserved quantities is illustrated
with the example of an axially symmetric top formulated as three
dimensional Nambu system.
\end{abstract}

\section{ Introduction }

Nambu introduced a generalization of Hamiltonian mechanics in
1973~\cite{Nambu}. In this generalization, points of the phase-space
were labeled by a canonical triplet $(x,y,z)$ and the evolution
equations were expressed in terms of a pair of Nambu functions $\Ho$
and $\Ht$, through
\begin{eqnarray}
  {d\vec{r} \over dt} = \nabla \Ho \times \nabla \Ht \nonumber 
\end{eqnarray}

In the subsequent years, the activities related to Nambu systems have
steadily developed~\cite{Mukunda,BAndF,KAndT,Takht,Esta,Fecko}. In
~\cite{Sag} we proposed a geometric formulation of generalized Nambu
systems. The motivation behind that work was to provide a framework
suitable from dynamical view point. Other generalizations, such
as~\cite{Takht}, are quite elegant. However, in them, the dynamics on
$n$ dimensional system, is governed by $(n-1)$ Nambu functions. i.e.,
too many integrals of motion are assumed. The framework proposed
in~\cite{Sag} involves a $3n$
dimensional manifold $\M$, together with a constant
(i.e. corresponding to a constant section in the bundle of 3-forms), strictly
non-degenerate  3-form $\omegath$. The time evolution is governed by a
pair (and not $(3n -1)$) of Nambu functions $\Ho, \Ht$. An interesting
feature of this 
formalism is the following. In addition to the three-bracket of functions, a
two bracket of 2-forms is defined. This bracket is called
Nambu bracket. The Nambu bracket of 2-forms gives rise to a {\em Lie
  algebra}, whereas, the three bracket of functions does not give rise 
to an algebra. It was found that the formulation involving 2-forms
provides a natural approach to a Nambu system. In the present paper we
show that the stringent condition of constancy of 3-form $\omegath$
can be moderated by a more general condition of closedness. Also we
argue that the formalism involving 2-forms, considered as dynamical
variables,  is the most suitable to express geometric facts about
Nambu systems. Equivalence classes of pairs of functions, which are
associated with 2-forms, can be chosen simply by fixing the gauge. The
phase point is determined (up to gauge) in terms of 2-forms considered
as dynamical variables.  We prove Nambu-Darboux theorem in this new
setup. we also discuss systematically the connection between
symmetries of Nambu systems and conserved 2-forms. In particular, we
introduce the notion of  ``momentum map'' and carry out the
generalization of symplectic Noether theorem for Nambu systems. Unlike
in the case of Hamiltonian systems, here the roles of generator and of 
conserved quantity are played by 2-forms. 

In section~\ref{NambuSystems} we list the main features of geometric
formulation developed in~\cite{Sag}, the new statement of
Theorem~\ref{ND_thm} incorporates more general systems as Nambu
systems. In section~\ref{MomentumMappings} we introduce momentum maps,
symmetries and Noether theorem for Nambu systems and illustrate them
with the help of an example of a symmetric top. 

\section{ Nambu Systems }
\label{NambuSystems}

In this section we will briefly review the relevant definitions and
results regarding Nambu systems. The Nambu structure was introduced through
the following definitions~\cite{Sag}. We also prove that even with
{\em closed} 3-form (instead of {\em constant}) a Nambu-Darboux
coordinates can be obtained.

\begin{defn} (Nondegenerate form) : Let $E$ be a finite dimensional
  vector space and let $\omegath$ be a 3-form on $E$ i.e., 
  \begin{eqnarray}
    \omegath : E \times E \times E \rightarrow \R  \nonumber
  \end{eqnarray}
  the form $\omegath$ is called a {\em nondegenerate form} if
  \begin{eqnarray}
    \forall\;\; non\;\;zero\;\;e_1 \in E\;\; \exists \;\; e_2, e_3 \in
    E\;\; \nonumber \\
    such \;\; that \;\; \omegath (e_1, e_2, e_3) \not= 0 \nonumber
  \end{eqnarray}
\end{defn}

\begin{defn} (Nambu complement) : Let $E$ be an $m$ dimensional vector
  space with $m \geq 3$. Let $\omegath$ be an anti-symmetric and
  non-degenerate 3-form on $E$. Let us choose $e_1,e_2,e_3 \in E$ such
  that $\omegath(e_1,e_2,e_3) \not= 0$. Let $P_1 = Span(e_1,e_2,e_3)$,
  then the {\em Nambu complement} of $P_1$ is defined as 
  \begin{eqnarray}
    P_1^{\perp_E} = \{ z \in E \;\;|\;\; \omegath(z,z_1,z_2) = 0
    \;\;\forall\;\;z_1,z_2 \in P_1 \} \nonumber 
  \end{eqnarray}
\end{defn}

\begin{defn} (Strictly nondegenerate form) : Let $E$ be an m
  dimensional vector space and $\omegath$ be an anti-symmetric and
  non-degenerate 3-form on $E$, the $\omegath$ is called
  {\em strictly non-degenerate} if for each non zero $e_1 \in E$
  $\exists$ a two dimensional subspace $E_1 \subset E$ such that  
  \begin{enumerate}
  \item $\omegath(e_1,x_1,x_2) \not= 0$ $\forall$ linearly independent
    $\{e_1,x_1,x_2\}$ where $x_1,x_2 \in F_1$ and $F_1 = Span(e_1 +
    E_1)$. 
  \item $\omegath(e_1,z_1,z_2)=0$ $\forall$ $z_1, z_2 \in
    F_1^{\perp_E}$. 
  \end{enumerate}
\end{defn}

Using the definition of Strictly nondegenerate 3-form, we define the Nambu
manifold as

\begin{defn} (Nambu Manifold) : Let $\M$ be a $3n$-dimensional $C^\infty$
  manifold and let $\omegath$ be a 3-form field on $\M$ such that
  $\omegath$ is completely anti-symmetric, closed and strictly nondegenerate at every
  point of $\M$ then the pair $(\M,\omegath)$ is called a
  {\em Nambu manifold}. 
\end{defn}
{\bf Remark:} Note the change in the definition from that
of~\cite{Sag}. The constant 3-form is now replaced by closed 3-form.
An analog of canonical transformation was defined in~\cite{Sag}.
\begin{defn} (Canonical transformation) : Let $(\M,\omegath)$ and
  $(N^{3n},\rho^\th)$ be Nambu manifolds. A $C^\infty$ mapping $F : \M
  \rightarrow N^{3n}$ is called {\em canonical transformation} if
  $F^*\rho^\th = \omegath$,
  where $F^*$ is the pullback of $F$.
\end{defn}

\begin{thm} (Nambu-Darboux theorem) : \label{ND_thm} Let
  $(\M,\omegath)$ be a Nambu manifold then  at every point $p \in \M$,
  there is a chart $(U,\phi)$ in which $\omegath$ is written as 
  \begin{eqnarray}
    \omegath |_U = \sum_{i=0}^{n-1} dx_{3i+1} \wedge dx_{3i+2} \wedge dx_{3i+3}
    \nonumber
  \end{eqnarray}
  where $(x_1,x_2,x_3,\ldots,x_{3(n-1)+1},x_{3(n-1)+2},x_{3(n-1)+3})$
  are local coordinates on $U$ described by $\phi$. 
\end{thm}
\pf
We wish to prove that there exists a canonical transformation which
maps any coordinate system to a system in which
$\omegath$ has the required form. Here we introduce the following
notation. Without loss of generality we prove the theorem on a vector
space $E$ with point $p$ as the {\em null vector} of the vector space
$E$. Let $\omega_1 = \omegath (0)$ be a constant 3-form on this
chart. We define a difference form $\tilde{\omega} = \omega_1 -
\omegath$. Let 
\begin{eqnarray}
  \omega_t &=& \omegath + t \tilde{\omega} \;\;\; 0 \leq t \leq 1
  \nonumber \\
  &=& (1-t) \omegath + t \omega_1 \nonumber
\end{eqnarray}
At $x=0$ we have $\tilde{\omega}=0$. Hence, $\omega_t (0) =
\omegath(0)$ is strictly non-degenerate $\forall$ $0 \leq t \leq
1$. Also for all $x$ the form $\omega_t$ is strictly non-degenerate at
$t=0$ and $t=1$.

To show that $\omega_t$ is strictly-non-degenerate for $0 < t < 1$.\\
Consider a three dimensional subspace $P$ of $E$. There are two cases \\
Case I: If $\omegath |_P = 0 \Rightarrow \omega_1|_P = 0 \Rightarrow
\omega_t|_P = 0$ \\ 
Case II: $\omegath|_P \not= 0 \Rightarrow \omega_1 |_P \not= 0
\Rightarrow \exists$ some neighborhood around $p=0$ such that
$\omega_t|_P \not= 0$ since $\omega_t$ is a continuous field.

Since $\omegath$ is strictly non-degenerate there are only two such
cases. So $\omega_t$ is also strictly non-degenerate.

${\omegath}$ is closed $\Rightarrow$ $\tilde{\omega}$ is closed.\\
Hence locally, by the Poincar\'e lemma, $\tilde{\omega} = d \alpha$ where
$\alpha$ is a 2-form so chosen that $\alpha(0) = 0$.

Since $\omega_t$ is strictly-non-degenerate we have $n$ three
dimensional subspaces of $E$ (say $P_i,\;\; i=1,\ldots,n$) on which
$\omega_t \not= 0$.
These subspaces are Nambu complements of each other. These are the
only subspaces on which $\omega_t \not= 0$. Hence, there exist vector
fields $X_{t_i}$ such that
\begin{eqnarray}
  i_{X_{t_i}} \omega_t |_{P_i} = - \alpha |_{P_i} \;\;\; \forall \;\;
  i=1,\dots,n \nonumber 
\end{eqnarray}
Since $E = P_1 \oplus \cdots \oplus P_n$\\
$X_t = \sum_{i=0}^{n-1} X_{t_i}$
which is zero at $x=0$. Let $F_t$ be the integral curve of $X_t$
\begin{eqnarray}
  {d \over dt} (F^*_t \omega_t) &=& F^*_t (L_{X_t} \omega_t) + F^*_t {d
    \over dt}\omega_t \nonumber \\ 
  &=& F^*_t d i_{X_t} \omega_t + F^*_t \tilde{\omega} \nonumber \\
  &=& F^*_t (-d\alpha + \tilde{\omega}) = 0 \nonumber
\end{eqnarray}
Hence, $F_1$ provides the required coordinate transformation \\
\qed

The theorem establishes the existence of a canonical coordinate system,
in which $\omegath$ has a ``Normal form''.

Let $\T^0_k (\M)$ denote a bundle of k-forms on $\M$, $\Omega^0_k
(\M)$ denote the space of k-form fields on $\M$ and $\X(\M)$ denote
the space of vector fields on $\M$. Now for a given vector field $X$
on $\M$, we define the inner product of $X$ with k-form (or contraction
of a k-form by $X$) as
\begin{eqnarray}
  (i_X \eta^{(k)}) (a_1,\ldots,a_{k-1}) = \eta^{(k)} (X,a_1, \ldots, a_{k-1})
  \nonumber
\end{eqnarray}
where $\eta^{(k)} \in \Omega^0_k(\M)$ and $a_1, \ldots, a_{k-1} \in \X(\M)$

The operations of Lowering and Raising were defined in~\cite{Sag} as
follows. The Lowering map $\flat : \X(\M) \rightarrow \Omega^0_2(\M)$
defined by $ X \mapsto X^\flat = i_X \omegath$, and the Raising
map $\sharp : \Omega^0_2(\M) \rightarrow \X(\M)$, is defined by the
following prescription.
Let $\alpha$ be a 2-form and $\alpha_{ij}$ be its components in
Nambu-Darboux coordinates, then the components of $\alpha^\sharp$ are
given by 
\begin{eqnarray} \label{2_form_to_vec}
  {\alpha^\sharp}^{3i+p} = {1 \over 2} \sum_{l,m=1}^3 \varepsilon_{plm}
  \alpha_{3i+l \;\;\; 3i+m}
\end{eqnarray}
where $0 \leq i \leq n-1$, $ p=1,2,3$ and $\varepsilon_{plm}$ is the
Levi-Cevit\`a symbol.

Here we have an important remark to make. On the space
$\T^0_{2_x}(\M)$, the space of 2-forms at $x \in \M$ the $\sharp$
defines an equivalence relation such as
$\omega_1^\tw(x) \sim \omega_2^\tw (x)$ if $(\omega_1^\tw)^\sharp (x)
= (\omega_2^\tw)^\sharp (x)$, where $\omega_1^\tw, \omega_2^\tw \in
\T^0_{2_x}(\M)$.

Equation (\ref{2_form_to_vec}) provides a relation
between 2-forms and vector fields. The following theorem, established
in~\cite{Sag}, provides conditions on 2-form, under which the
associated flow preserves the Nambu structure.
\begin{thm}\label{equiv_close_form}
  Let $\beta^{(2)} \in \Omega^0_2(\M)$, and $f^t$ be a flow
  corresponding to $\beta^{(2)^\sharp}$, 
  i.e., $f^t : \M \rightarrow \M$ such that 
  \begin{eqnarray}
    {\div d { dt } } \Big|_{t=0} ( f^tx ) = (\beta^{(2)^\sharp}) x \;\;
    \forall x \in \M
    \nonumber 
  \end{eqnarray}
  Then the form $\omegath$ is preserved under the action of
  $\beta^{(2)^\sharp}$ iff $d(\beta^{(2)^\sharp})^\flat = 0$.
  i.e., ${f^t}^* \omegath = \omegath$ iff $d(\beta^{(2)^\sharp})^\flat = 0$ 
\end{thm}

From theorem~\ref{equiv_close_form} it follows that the vector fields
preserving the Nambu structure can be obtained from two functions,
called Nambu functions, as in~\cite{Sag}. These Nambu functions were shown to
be constants of motion. 
At this stage, the following definitions were introduced.
\begin{defn} (Nambu vector field) : Let ${\Ha 1 }$, ${\Ha 2 }$ be real valued 
  $C^\infty$ functions (Nambu functions) on $(\M,\omegath)$ then $N$ is
  called {\em Nambu vector field} corresponding to ${\Ha 1 }, {\Ha
    2 }$ if 
  \begin{eqnarray}
    N = (d {\Ha 1 } \wedge d {\Ha 2 })^\sharp \nonumber
  \end{eqnarray}
\end{defn}

\begin{defn} (Nambu system) : A four tuple $(\M,\omegath, {\Ha 1 },
  {\Ha 2 })$ is called {\em Nambu system}. 
\end{defn}

\begin{defn} (Nambu phase flow) : Let $(\M, \omegath , {\Ha 1 }, {\Ha 2 })$ 
  be a Nambu system and let the $g^t : \M \rightarrow \M$ be one
  parameter family of group of diffeomorphisms. The flow corresponding
  to $g^t$ is
  {\em Nambu flow} if 
  \begin{eqnarray}
    {\div d {dt} } \Big|_{t=0} (g^t {\bf x}) &=& (d{\Ha 1 } \wedge d {\Ha
      2 })^\sharp {\bf x} \;\;\; \forall\; {\bf x} \in \M \nonumber \\ 
    &=& N {\bf x} \nonumber
  \end{eqnarray}
\end{defn} 

It can be shown that $g^t$ preserves $\omegath$~\cite{Sag}.

\begin{defn} (Nambu bracket) : \label{NambuBracket} Let $\omega_a^\tw$
  and $\omega_b^\tw$ be 2-forms then the {\em Nambu bracket} is a
  map $\{,\} : \Omega^0_2(\M) \times \Omega^0_2(\M) \rightarrow
  \Omega^0_2(\M)$ given by 
  \begin{eqnarray}
    \{ \omega_a^\tw , \omega_b^\tw \} = [ {\omega_a^\tw}^\sharp ,
    {\omega_b^\tw}^\sharp ]^\flat \nonumber 
  \end{eqnarray}
  where $[,]$ is Lie bracket of vector fields.
\end{defn}

From definition~\ref{NambuBracket} it follows that the bracket of
2-forms provides a Lie algebra of 2-forms.

In case of Hamiltonian systems, the Poisson bracket generates the
Hamiltonian vector field. A similar result was proved in~\cite{Sag}.
This assures that if $\a$ and $\b$ are 2-forms, if 
$\a^\sharp$ is  a Nambu vector field and if $\a^\prime = (\a^\sharp)^\flat$ 
and $\b^\prime = (\b^\sharp)^\flat$ then 
\begin{eqnarray}
  \{ \a, \b \} = L_{{\a^\prime}^\sharp} {\b^\prime} \nonumber
\end{eqnarray}
where  $L_{{\a^\prime}^\sharp}$ denote the Lie derivative along the
vector field ${\a^\prime}^\sharp$.

\begin{defn} (Nambu bracket for functions) : Consider a Nambu manifold 
  $(\M, \omegath)$ and let $f,g,h$ be $C^{\infty}$ functions on $\M$
  then {\em Nambu bracket for functions} is given by 
  \begin{eqnarray}
    \{f,g,h \} = L_{(dg \wedge dh)^\sharp}f = i_{(dg \wedge dh)^\sharp} df
    \nonumber 
  \end{eqnarray}
\end{defn}

We wish to emphasize that in case of the Hamiltonian systems, the bracket
of functions is equivalent to  the bracket of 1-forms. In the present
formulation, the bracket of 2-forms and the bracket of functions are not
equivalent, but they are related to each other as given in the
following proposition which is proved in~\cite{Sag}.

\begin{prop}\label{bracket_relation}
  Let $(\M,\omegath)$ be a Nambu manifold and let $f,g,h_1,h_2$ be
  $C^\infty$ functions satisfying $(df \wedge dg)^{\sharp^\flat} = df
  \wedge dg$ and $(dh_1 \wedge dh_2)^{\sharp^\flat} = dh_1 \wedge dh_2$
  then 
  \begin{eqnarray}
    \{ dh_1 \wedge dh_2, df \wedge dg \} = d \{f,h_1,h_2\} \wedge dg + df
    \wedge d\{g,h_1,h_2\} \nonumber 
  \end{eqnarray}
\end{prop}

\section{ Nambu momentum maps }
\label{MomentumMappings}

\subsection {Forms as dynamical variables}

Let us recall the Hamiltonian formalism. In that case, the phase
space is a symplectic manifold $(M^{2n},\omegatw)$. i.e., an even
dimensional manifold equipped with a closed non-degenerate
2-form. The vector field are in one to one correspondence with
the 1-forms as $\a = \omegatw(\cdot,X_\a)$ where $\a \in
\Omega^0_2(M^{2n})$ and $X_\a \in \X(M^{2n})$~\cite{AM}. Hence,
Hamiltonian systems can be completely formulated using 
1-forms considered as dynamical variables. Consequently the
conservation laws and symmetry are also formulated in terms of
1-form. The relation between the 1-form and functions involve a one
parameter family of freedom (e.g.: adding a constant to Hamiltonian
does not change the equations of motion). Thus, for convenience, the
constants of motion are considered as functions rather than 1-forms.

Now, in the Nambu framework, it is natural to consider 2-form fields
as dynamical variables. Because, as mentioned above, there is a
correspondence between vector fields and classes of 2-forms~\cite{Sag}.
In fact, it is possible to obtain other types of dynamical variables
from 2-forms using the following procedure. As shown in~\cite{Sag},
the space of equivalence classes of 2-forms at each point is a 3n-dimensional
space. A set of $3n$ independent vector fields, and subsequently $3n$ local
coordinates can be obtained from $3n$ independent 2-forms. Thus all types of
dynamical variables (e.g. functions) can be recovered, in principle, from
the convenient set of 2-forms. 

We call a pair of functions $(f,g)$ to be equivalent to another pair
$(h,p)$ if $df \wedge dg = dh \wedge dp$. So the equivalence classes
of pairs of functions are locally associated with closed 2-forms. Any
closed 2-form can be represented by giving a representative pair in
the equivalence classes of pairs of functions. The choice of
representative pair is called ``fixing the gauge'' (See~\cite{Sag}).

\subsection{Nambu momentum map}
Before formalizing the concept of Nambu momentum mapping we state the
standard definitions.

\begin{defn} (Invariants~\cite{AM}) : Let $\M$ be a manifold and $X\in
  \X(\M)$. Let $\a$ be a k-form on $\M$. We call $\a$ an invariant
  k-form with respect to $X$ iff $L_X \a = 0$. 
\end{defn}

\begin{defn} (Group action on manifold~\cite{AM}) : Let $M$ be a
  manifold and let $G$ be a Lie group. A {\em group action on
    manifold} $M$ is a smooth mapping $\Phi_G : G \times M \rightarrow
  M$ such that  
  \begin{enumerate}
  \item $\Phi_G ( e, x) = x$ 
  \item $\Phi_G ( g, \Phi_G (h,x)) = \Phi_G(gh,x)$ $\forall$ $g,h \in
    G$ and $ x \in M$ 
  \end{enumerate}
\end{defn}

\begin{defn} (Infinitesimal generator~\cite{AM}) : Let $\Phi_G : G
  \times M \rightarrow M$ be action of group $G$ on manifold $M$. For
  $\xi \in \g$ where $\g$ is Lie algebra of $G$, the map $\Phi_G^\xi :
  \R \times M \rightarrow M$, defined by $\Phi_G^\xi(t,x) =
  \Phi_G(\exp (t\xi), x)$ is an $(\R,+)$ action on $M$. The {\em
    infinitesimal generator}, $\xi_M \in \X(M)$, of action
  corresponding to $\xi$ is  
  \begin{eqnarray}
    \xi_M (x) = {d \over dt} \Phi_{\exp(t\xi)} (x) \Big|_{t=0} \nonumber
  \end{eqnarray}
\end{defn}

Before proceeding further with Nambu-Noether theorem we briefly
recall the corresponding ideas in Hamiltonian dynamics~\cite{AM}. A
group action on symplectic manifold is called a symplectic action if
it preserves the symplectic structure. A $\g^*$ valued map, (where
$\g^*$ is dual of Lie algebra), on symplectic manifold is called a
momentum map, provided it satisfies certain consistency conditions as
shown in~\cite{AM}.

We now develop similar ideas for Nambu systems.
\begin{defn} (Nambu action) : Let $G$ be a Lie group and let
  $(\M,\omegath)$ be a Nambu manifold. Let $\Phi_G$ be the action of
  $G$ on $\M$. $\Phi_G$ is called {\em Nambu action} if 
  \begin{eqnarray}
    \Phi_G^* \omegath = \omegath \nonumber
  \end{eqnarray}
  i.e., $\Phi_G$ induces a Nambu canonical transformation.
\end{defn}

We introduce a 2-form
valued quantity called momentum, which is a generator of the group
action, associated with each one-parameter group of symmetry. Here we
define the ``Nambu momentum map'': 
\begin{defn} (Nambu momentum maps) : Let $G$ be a Lie group and let
  $(\M, \omegath)$ be a Nambu manifold, let $\Phi_G$ be a Nambu action on
  $\M$. Then the mapping ${\bf J} \equiv (J_1, J_2) : \M \rightarrow
  \g^* \times \g^*$ is called {\em Nambu momentum maps} provided
  for every $\xi \in \g$ 
  \begin{eqnarray}
    (d\hat{J}_1 (\xi) \wedge d\hat{J}_2 (\xi))^\sharp = \xi_{\M} \nonumber
  \end{eqnarray}
  where $\hat{J}_1(\xi) : \M \rightarrow \R$, $\hat{J}_2 (\xi) :
  \M\rightarrow\R$, defined by $\hat{J}_1 (\xi)(x) = J_1(x) \xi$ and
  $\hat{J}_2(\xi)(x)= J_2(x)\xi$ $\forall x \in\M$ and $\xi_{\M}$ is
  infinitesimal generator of the action $\Phi_G$.
\end{defn} 
\noindent
{\bf Remark: } 
\begin{enumerate}
\item Traditionally the term momentum map has come from the conjugate
momenta. Here we can not associate any such meaning to this quantity. We
call it as Nambu momentum map because under symetry the 2-form constructed
from Nambu momentum map is preserved under Nambu flow (see
Theorem~\ref{Noether}) like the usual momentum.
\end{enumerate}
\begin{defn} (Nambu G-space) : The five tuple
  $(\M,\omegath,\Phi_G,J_1,J_2)$ is called {\em Nambu G-space}. 
\end{defn}

The following proposition establishes the consistency between the
bracket of 2-forms and the Nambu momentum maps.
\begin{prop}\label{Consistency} Let $(\M,\omegath,\Phi_G, J_1, J_2)$ be
  a Nambu G-space and $\xi,\eta \in \g$ then 
  \begin{eqnarray}
    (d\hat{J}_1([\xi,\eta]) \wedge d\hat{J}_2([\xi,\eta]))^\sharp = \{
    d\hat{J}_1(\eta) \wedge d\hat{J}_2(\eta), d\hat{J}_1(\xi) \wedge
    d\hat{J}_2(\xi) \}^\sharp 
    \nonumber
  \end{eqnarray}
i.e., The following diagram commutes \\
\begin{center}
\begin{picture}(100,100)
\put(5,80){$\Omega^0_2(\M)$}
\put(145,80){$\X(\M)$}
\put(5,10){$\g$}
\put(70,87){$\sharp$}
\put(-17,40){$\hat{J_1}, \hat{J_2}$}
\put(72,40){$\xi \mapsto \xi_{\M}$}
\thicklines
\put(50,82){\vector(1,0){90}}
\put(9,17){\vector(0,1){58}}
\put(11,17){\vector(2,1){129.4}}
\end{picture}
\end{center}
\end{prop}
\pf 
\begin{eqnarray}
  \{ d\hat{J}_1(\eta) \wedge d\hat{J}_2(\eta)&,& d\hat{J}_1(\xi) \wedge
  d\hat{J}_2(\xi) \}^\sharp \nonumber \\
  &=& [ (d\hat{J}_1(\eta) \wedge d\hat{J}_2(\eta))^\sharp,
  (d\hat{J}_1(\xi) \wedge d\hat{J}_2(\xi))^\sharp ] \nonumber \\ 
  &=& [\eta_{\M}, \xi_{\M} ] \nonumber \\
  &=& - [\eta, \xi]_{\M} \nonumber \\
  &=& (d\hat{J}_1([\xi,\eta]) \wedge d\hat{J}_2([\xi,\eta]))^\sharp \nonumber
\end{eqnarray}
\\ \qed 

\noindent
{\bf Remark:}
\begin{enumerate}
\item Infact, if one starts with bracket of 2-form as more basic
  quantity than the closed, strictly non-degenerate 3-form $\omegath$, 
  then the Proposition~\ref{Consistency} can be used as the definition 
  of the Nambu momentum maps. The equivalence of such framework with
  that of current one is not yet established.
\end{enumerate}
In the following examples, we explicitly construct Nambu momentum maps
for the group actions of $SO(3)$ and $SP(2)$. They also illustrate the
consistency condition stated in Proposition~\ref{Consistency}\\
\ex Consider action of $SO(3)$ on the Nambu manifold $(\R^3, dx \wedge
dy \wedge dz)$, where $(x,y,z)$ are the Nambu-Darboux coordinates on
$\R^3$. The Lie algebra $\g = \R^3$. Let $e_1, e_2, e_3$ be a basis of 
$\g$ satisfying the bracket conditions. 
\begin{eqnarray} 
  [ e_1, e_2 ] = e_3,\;\;
  [ e_2, e_3 ] = e_1,\;\;
  [ e_3, e_1 ] = e_2 \nonumber
\end{eqnarray}
Let $f_1, f_2, f_3$ be a dual basis of $\g^*$ corresponding to $e_1, e_2, e_3$.

Let $\vec{r} \equiv x f_1 + y f_2 + z f_3, \vec{\rho} \equiv {1 \over 2} (y^2 +
z^2) f_1 + {1 \over 2} (x^2 + z^2) f_2 + {1 \over 2} (y^2 + x^2) f_3 \in \g^*$ then
the Nambu momentum maps in the Nambu-Darboux coordinates are
$J_1(x,y,z) = \vec{r}$ and $J_2(x,y,z) = \vec{\rho}$. The 2-forms
obtained from these are  
\begin{eqnarray}
  d\hat{J}_1(e_1) \wedge d\hat{J}_2(e_1) &\equiv& L_1 = y dx \wedge dy + z dx
  \wedge dz \nonumber \\ 
  d\hat{J}_1(e_2) \wedge d\hat{J}_2(e_2) &\equiv& L_2 = - x dx \wedge dy + z
  dy \wedge dz \nonumber \\ 
  d\hat{J}_1(e_3) \wedge d\hat{J}_2(e_3) &\equiv& L_3 = - x dx \wedge dz - y
  dy \wedge dz \nonumber 
\end{eqnarray}

Little algebra verifies that $\{ L_1, L_2 \} = - L_3$, $\{ L_2, L_3\}
= - L_1$
and $\{L_3, L_1 \} = - L_2$. This establishes the consistency with
Proposition~\ref{Consistency}. \\
\ex Consider the action of $SP(2,\R)$ on $(\R^3, dx \wedge dy \wedge
dz)$ defined by $x \mapsto x$ and $(y,z) \mapsto A \cdot (y,z)$ where
$A \in SP(2,\R)$ and $(x,y,z)$ are the Nambu-Darboux coordinates on
$\R^3$. It is easy to check that the action of $SP(2,\R)$ is 
a Nambu action\footnote{ It follows from the fact that the $SP(2,\R)$
preserves area on the two dimensional plane $(y,z)$ and the action is
identity on the orthogonal $x$ direction. Hence, the three volume is
preserved under this action.}. The algebra $sp(2,\R)$ is three
dimensional. Let  
\begin{eqnarray}
  e_1 = \left[ \begin{array}{cc}
      1 & 0 \\
      0 & -1
    \end{array} \right],
  e_2 = \left[ \begin{array}{cc}
      0 & -1 \\
      1 & 0
    \end{array} \right],
  e_3 = \left[ \begin{array}{cc}
      0 & 1 \\
      1 & 0
    \end{array} \right] \nonumber
\end{eqnarray}
be a basis of $sp(2,\R)$. The bracket relations are
\begin{eqnarray}
  [ e_1, e_2] = -2 e_3, [e_2, e_3 ] = -2 e_1, [e_3, e_1] = 2 e_2 \nonumber
\end{eqnarray}
Let $f_1, f_2, f_3$ be a dual basis of $\g^*$ corresponding to $e_1, e_2, e_3$.

The Nambu momentum maps in Nambu-Darboux coordinates are
\begin{eqnarray}
  J_1(x,y,z) &=& z f_1 + y f_2 + x f_3 \in sp^*(2,\R) \nonumber \\
  J_2(x,y,z) &=& (x^2 - y^2) f_1 +  (x^2 + z^2) f_2 +  (y^2 - z^2) f_3\in sp^*(2,\R_) \nonumber
\end{eqnarray}

Let $A_i = d\hat{J}_1(e_i) \wedge d\hat{J}_2(e_i), \;\; i=1,2,3$. Thus 
\begin{eqnarray}
  A_1 &=& -2x dx \wedge dz + 2y dy \wedge dz \nonumber \\
  A_2 &=& -2x dx \wedge dy + 2z dy \wedge dz \nonumber \\
  A_3 &=& 2y dx \wedge dy - 2z dx \wedge dz \nonumber
\end{eqnarray}

It is easy to check that $A_1, A_2, A_3$ satisfy the consistency
condition established in proposition~\ref{Consistency}. It
so happens that the action of $SP(2)$ just displaces the points of
$\R^3$ on some conic section.

In analogy with the notion of symplectic symmetry, we now introduce
Nambu Lie symmetry. 
\begin{defn} (Nambu Lie symmetry) : Consider a Nambu system
  $(\M,\omegath,\Ho,\Ht)$. Let 
  $(\M,\omegath,\Phi_G,J_1,J_2)$ be a Nambu G-space. We call $\Phi_G$ a
  {\em Nambu Lie symmetry} of the Nambu system if
  \begin{eqnarray}
    \Phi_G^* (d\Ho \wedge d\Ht) = (d\Ho \wedge d\Ht) \nonumber
  \end{eqnarray}
\end{defn}

We now come to the main theme of this paper. We first prove
the theorem which establishes a connection between the symmetries and
invariants of Nambu systems, i.e., Noether theorem for Nambu systems
and then elaborate the idea of the theorem, by giving an example of
Symmetric top. 
\begin{thm} (Nambu-Noether theorem): \label{Noether}Consider a Nambu
  system $($ $\M, \omegath, \Ho, \Ht$ $)$ where 
  $\Ho, \Ht$ are so chosen that $d\Ho \wedge d\Ht = {(d\Ho \wedge
    d\Ht)^\sharp}^\flat$. Let this system be a Nambu G-space
  $(\M,\omegath,\Phi_G,J_1,J_2)$  where $J_1, J_2$ are so chosen that
  $d\hat{J}_1(\xi) \wedge d\hat{J}_2(\xi) = {(d\hat{J}_1(\xi) \wedge
    d\hat{J}_2(\xi))^\sharp}^\flat$ $\forall \xi \in \g$. 
  If $\Phi_G$ is Nambu Lie symmetry of this system
  then
  \begin{eqnarray}
    L_{(d\Ho \wedge d\Ht)^\sharp} \Big(d\hat{J}_1(\xi) \wedge d\hat{J}_2(\xi)\Big)
    = 0 \nonumber 
  \end{eqnarray}
  i.e., $d\hat{J}_1(\xi) \wedge d\hat{J}_2(\xi)$ is conserved by  the
  Nambu flow.
\end{thm}
\pf Let $g^t = \exp(t\xi)$ be a one parameter subgroup of $G$. In view
of the choice of $\Ho$, $\Ht$ and Nambu Lie symmetry it follows that 
\begin{eqnarray}
  0 &=& {d \over dt} (\Phi_{\exp(t\xi)}^* (d\Ho \wedge d\Ht)) \Big|_{t=0} \nonumber \\ 
  &=&\Phi_{e}^* (L_{\xi_{\M}} (d\Ho \wedge d\Ht)) \nonumber \\
  &=&\Phi_{e}^*\{ \xi_{\M}^\flat, d\Ho \wedge d\Ht \} \nonumber
\end{eqnarray}
Hence
\begin{eqnarray}
  0 &=& \Phi_{e}^*\{  d\Ho \wedge d\Ht, \xi_{\M}^\flat \} \nonumber \\
  &=& L_{(d\Ho \wedge d\Ht)^\sharp} (d\hat{J}_1(\xi) \wedge d\hat{J}_2(\xi))
  \nonumber 
\end{eqnarray}
which follow from the definition of Nambu momentum maps and the choice of
$\hat{J}_1$ and $\hat{J}_2$.\\
\qed \\
{\bf Remark:} \\

In Hamiltonian dynamics the vector field are in one to one
correspondence with the 1-forms as $\a = \omegatw(\cdot,X_\a)$ where
$\a \in \Omega^0_2(M^{2n})$ and $X_\a \in
\X(M^{2n})$~\cite{AM}. Since 1-forms are locally related to the functions, 
up to a constant, it is traditional to consider functions as generators
and conserved quantities. However, 1-forms play role of generators.  
In the present framework the roles of generator and of conserved
quantity is played by 2-form.\\
\ex Consider the action of $SO(2)$ on the manifold $(\R^3, dx\wedge dy
\wedge dz)$ defined by $x \mapsto x$, $(y,z) \mapsto (y \cos(\phi) + z
\sin(\phi), -y \sin(\phi) + z \cos(\phi))$ where $\phi$ is angle of
rotation around $x$ axis and $(x,y,z)$ are the Nambu-Darboux
coordinates on $\R^3$. The Lie algebra of $SO(2)$ is $so(2) =
\R$. Let $e_1$ be the basis of $so(2)$. Let $f_1$ be a dual basis
corresponding to $e_1$. The Nambu momentum maps
corresponding to the action of $SO(2)$ in Nambu-Darboux coordinates
are
\begin{eqnarray}
  J_1(x,y,z) &=& x f_1 \nonumber \\
  J_2(x,y,z) &=& ({1 \over 2} ( y^2 + z^2 )) f_1\nonumber
\end{eqnarray}
The the 2-form 
\begin{eqnarray}
  d\hat{J}_1(e_1) \wedge d\hat{J}_2(e_1) \equiv  L_1 = y dx \wedge dy + z dx
  \wedge dz \nonumber 
\end{eqnarray}

Now consider a Nambu system for a symmetric top. Let $(x,y,z)$ be
the components of angular momentum of the top in the body frame. The
Nambu functions are 
\begin{eqnarray}
  \Ho &=& {1 \over 2} ( {x^2 \over I_x^2} + {y^2 \over I_y^2} + {z^2
    \over I_y^2}) \nonumber \\ 
  \Ht &=& x^2 + y^2 + z^2 \nonumber
\end{eqnarray}
It follows that $L_{L_1^\sharp} (d\Ho \wedge d\Ht) = 0$.
It is easy to see that $L_{(d\Ho \wedge d\Ht)^\sharp} L_1 = 0$. The
$L_1$ is conserved.

\section{ Conclusions }

In the present paper we have enlarged the class of Nambu systems. This
is achieved by replacing the stringent condition of constancy of 3-form 
by the closedness of the 3-form. We have also developed the notion of
symmetry of a Nambu system. In the geometric formulation of Nambu
mechanics, it is natural to consider 2-forms as dynamical variables. 
Therefore, the notion of momentum maps
in this framework is formulated in terms of 2-forms. The
connection between continuous groups of symmetries and the
conservation laws, i.e., the Noether theorem, is proved in the present
work. The result is  illustrated with the help of an example of 
an axially symmetric 3-dimensional Nambu system. 

\ack
We thank Dr. H. Bhate and Prof. K. B. Marathe for the critical comments
on the manuscript. One of the authors (SAP) is grateful to CSIR
(India) for financial assistance during the initial stages of work. 

\end{document}